\documentclass[11pt,a4paper]{article}
\usepackage{authblk}
\usepackage{amsmath}
\usepackage{amsfonts} 
\usepackage{hyperref}
\usepackage{graphicx}
\usepackage{subfigure}     % required only when side-by-side / subfigures are used

\newtheorem{problem}{Problem}
\newtheorem{definition}{Definition}

\begin{document}

\title{Topological data analysis and clustering}

\author[*]{Dimitrios Panagopoulos}

\affil[*]{Pelopa 3,\,Gerakas Attikis
Greece 15344,
dpanagop [at] gmail.com}
\maketitle

\begin{abstract}
Clustering is one of the most common tasks of Machine Learning. In this chapter we examine how ideas from topology can be used to improve clustering techniques.
\end{abstract}

%\tableofcontents

\section{Introduction}\label{panagopoulos_sec1}
With the advent of Big Data, algorithms that try to extract information from them are ubiquitous. Clustering algorithms are a subcategory of Machine Learning algorithms with a wide range of applications. Notions like closeness, distance and shape are central to clustering. It is then natural to try to use ideas and techniques from topology to improve clustering algorithms. This chapter examines some ways on how this could happen. In Section \ref{panagopoulos_sec2} a brief introduction to the clustering task is presented. In Subsection \ref{panagopoulos_subsec2_1} a definition is presented along with notation. k-means algorithm will be one of our focus areas and is examined in more detail in Subsection \ref{panagopoulos_subsec2_2}. In Section \ref{panagopoulos_sec3} Topological Data Analysis is presented. We will be using Persistent Homology and the Mapper algorithm which are presented in Subsections \ref{panagopoulos_subsec3_1} and \ref{panagopoulos_subsec3_2} respectively. Ways of applying topology to clustering are examined in Section \ref{panagopoulos_sec4}. Finally, Section \ref{panagopoulos_sec5} has some concluding remarks.

All figures in this article were created using Scikit-TDA package \cite{ST} in Python. The relative code can be found in \url{https://github.com/dpanagop/data_analytics_examples/tree/master/topological_data_analysis}.

\section{Clustering}\label{panagopoulos_sec2}
\subsection{Definition}\label{panagopoulos_subsec2_1}
Clustering is the task of grouping a set of objects in such a way that objects in the same group (called a cluster) are more similar (in some sense) to each other than to those in other groups \cite{Wiki}. Clustering is a major task of Machine Learning with diverse applications from medicine \cite{Sot2003} to sports \cite{Akhanli}.

The definition is intentionally general and is very common when one tries to apply clustering algorithms a lot of thought to be devoted on how similarity should be defined. For example, in marketing a common application is customer clustering. There one should consider if similarity depends on things like age, gender, place of residence, education e.t.c. It must be clear that in most cases there is not a predetermined grouping of the objects to be clustered against which an algorithm can be measured\footnote{In Machine Learning, the category of problems where there is no given labelling or values that can be used to train a model is called unsupervised learning.}.

In order to have some notation available we state the following problem:
\begin{problem}[Clustering Problem]\label{clustering_problem}
Given a finite subset $X\subset\mathbb{R}^d$ find a partition $C_1,\ldots,C_k$ of $X$ such that elements in every subset $C_i$ are more similar to other elements in $C_i$ than to elements in another subset $C_j$, where $j\neq i$. 
\end{problem}

Some of the most common algorithms used for clustering are:
\begin{enumerate}
\item k-means,
\item hierarchical clustering,
\item density-based spatial clustering of applications with noise (DBSCAN) and
\item Gaussian Mixture Models.
\end{enumerate}
In the next section k-means will be described in more detail. For the rest, the interested reader can start by reading the related Wikipedia article \cite{Wiki}, the documentation of Python's Scikit-learn package \cite{Pedregosa} at \url{https://scikit-learn.org/stable/modules/clustering.html} and \cite{YMJ}.

\subsection{K-means}\label{panagopoulos_subsec2_2}
One of the most popular algorithms for clustering is k-means. It is an algorithm that finds an approximate solution to the next problem.
\begin{problem}[k-means Problem]\label{k-means_problem}
Given a finite subset $X\subset\mathbb{R}^d$ and a natural number $k>1$, find a partition $C_1,\ldots,C_k$ of $X$ that minimizes
\begin{equation}
\sum_{i=1,\ldots,k}\sum_{x\in C_i}||x-c_i||^2,
\end{equation}\label{k-means_eq1}
where $c_i=\frac{1}{|C_i|}\sum_{x\in C_i}x$ is the barycenter of $C_i$.
\end{problem}

While k-means problem is known to be NP-hard \cite{ADHP,MNV}, the k-means algorithm is an efficient heuristic that detects a local minimum. The algorithm starts with a random selection of $k$ centers $c_1,\dots,c_k$. Then:
\begin{enumerate}
\item each $x\in X$ is assigned to its closest center, thus 
\[C_i=\{x: ||x-c_i||<||x-c_j||,\,j\in\{1,\ldots,k\}\setminus\{i\}\},\,i=1,\ldots,k\]
\item new centers are calculated based on equation $c_i=\frac{1}{|C_i|}\sum_{x\in C_i}x,\,i=1,\ldots,k$
\item steps 1 and 2 are repeated until centers $c_1,\dots,c_k$ do not change or a maximum number of iterations is reached. 
\end{enumerate}

\begin{figure}
\centerline{\includegraphics[width=6cm]{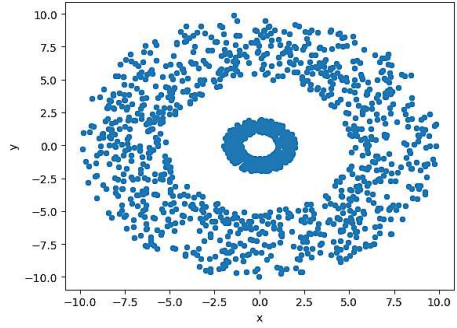}}
\caption{The image shows 1500 points in the plane. The points were selected at random such that 500 lie in a ring with inner radius 1 and outer 2. The  rest 1000 were selected such that they lie a ring with inner radius 5 and outer 10.}
\label{two_circles}
\end{figure}

There are several variations to avoid reaching a local minima and/or helping the algorithm converge faster. Those variations consist mainly of changing the way the initial centers are selected (ex. see \cite{AV,ORSS}). Furthermore, there are extensions for dealing with datasets $X$ with great number of points \cite{Sculley} or for replacing centers with real data points \cite{KR}.

k-means is a very efficient algorithm that converges fast and has, relatively, low memory requirements when implemented in a computer system. Two drawbacks of k-means are:
\begin{itemize}
\item the number of clusters (value of $k$) cannot be determined by the algorithm, instead it is an input for the algorithm,
\item k-means can detect clusters that have a spherical shape\footnote{It can be shown that there is no algorithm that is unaffected from scaling and from the act of packing more densely elements of the same cluster that can detect all possible partitions (i.e. that can detect clusters of arbitrary shapes) \cite{Kleinberg}.}.
\end{itemize}
For example, in fig.\ref{two_circles} a natural way to cluster the points would be in two concentric rings. k-means cannot achieve this. It is also worth mentioning, that k-means is not the only clustering algorithm that requires as input the number of clusters.  

\section{Topological Data Analysis}\label{panagopoulos_sec3}
Assuming that $X$ is a sample from a manifold, the number clusters is the number of connected components of the manifold. Furthermore, the shape of the connected components is related to the topological features of the manifold. Thus, it is natural to turn to the mathematical field of Topology for ideas to deal with the above problems. This gives rise to the field of Topological Data Analysis (TDA). In this section we present two of its tools, Topological Persistence and the Mapper algorithm. 

\subsection{Persistent Homology}\label{panagopoulos_subsec3_1}
Topological Persistence was introduced by Edelsbrunner, Letscher and Zomorodian in 2002 \cite{ELZ}. Let $X=\{x_1,\ldots,x_n\}\subset\mathbb{R}^d$ be a finite set of points. One can use $X$ to create a simplicial complex as below. 

\begin{definition}[Vietoris-Rips complex]\label{Vietoris_Rips}
Let $X=\{x_1,\ldots,x_n\}\subset\mathbb{R}^d$ be a finite set of points and $\epsilon>0$ a positive real number. The Vietoris-Rips complex $\mathcal{C}_{VR}(X,\epsilon)$ is the complex with $X$ as vertex set and where the n-simplex $\{x_{i_0},\ldots,x_{i_n}\}$ belongs to $\mathcal{C}_{VR}(X,\epsilon)$ if and only if $||x_{i_s}-x_{i_t}||<\epsilon$ for all $0\leq s,t\leq n$. 
\end{definition}

The Vietoris-Rips complex can be seen as an approximation of the manifold $X$ is sampled from and the parameter $\epsilon$ as the maximum distance between two points that are supposed to be path connected. 
 
\begin{figure}
\centerline{\includegraphics[width=8cm]{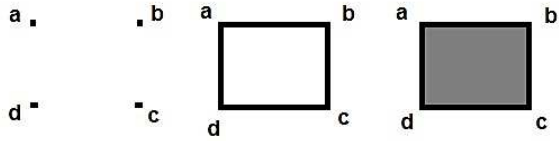}}
\caption{Four points on the corners of a unit square and the corresponding Vietoris-Rips complex for various selections of $\epsilon$. To the left when $\epsilon\in (0,1)$. In the middle when $\epsilon\in[1,\sqrt{2})$. And to the right when $\epsilon\geq\sqrt{2}$. }
\label{rips}
\end{figure}

A standard construction of Algebraic Topology, takes a simplicial complex $\mathcal{C}$ and constructs a family of abelian groups and connecting homomorphisms 
\[\cdots\rightarrow C_{n+1}\xrightarrow{\partial_{n+1}} C_n \xrightarrow{\partial_n} C_{n-1}\rightarrow\cdots\rightarrow C_1 \xrightarrow{\partial_1} C_0\xrightarrow{\partial_0}0\]
with $\partial_n\partial_{n+1}=0$ i.e $Im\partial_{n+1}\subseteq Ker\partial_n$ and where $C_n$ has as a base the n-simplices of $\mathcal{C}$. Thus, one can define the {\bf n-th homology group}\label{n_th_homology_group}
\[H_n:=Ker\partial_n/Im\partial_{n+1}\]
It can be shown if a complex has a finite number of n-simplices then the n-th homology group is a finitely generated abelian group. Thus, in this case the n-th homology group is isomorphic to the direct sum of copies of $\mathbb{Z}$ and finite cyclic groups, i.e. $H_n\simeq \mathbb{Z}^b\bigoplus_{i=1,\ldots,m} \mathbb{Z}_{p_i}$. The number of copies $\mathbb{Z}$ appears (i.e. $b$) is called the {\bf n-th betti number} $\beta_n$\label{n_th_betti_number}. Intuitively, $b_0$ counts the number of connected components, $b_1$ the number of 1-dimensional holes, $b_2$ the number of 2-dimensional holes e.t.c. For a more detailed exposition the reader can see \cite[Ch. 2]{Hatcher}.

\begin{figure}
\centerline{\includegraphics[width=6cm]{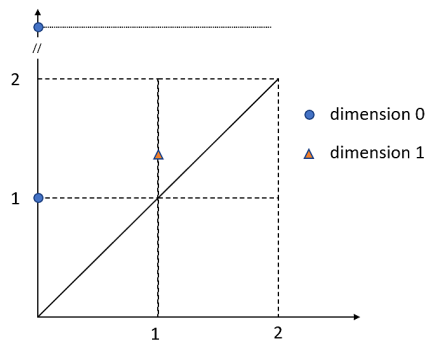}}
\caption{Persistent diagram for the Vietoris-Rips complexes of fig.\ref{rips}. There are 3 points at the $(0,1)$ that represent three of the four connected components when $\epsilon<1$. Since for $\epsilon=1$ all four vertices are connected, three of the four generators of $H_1$ vanish. The remaining generator does not vanish. This corresponds to the point that lies in the intersection of $x$-axis and the horizontal line at the top that represents infinity. For dimension 1 there is a loop that is created when $\epsilon=1$ and destroyed when $\epsilon=\sqrt{2}$. This is corresponds to the point that is marked with a triangle.}
\label{rips_diagram}
\end{figure}

It is clear that if $\epsilon<\epsilon'$ then every n-simplex of $\mathcal{C}_{VR}(X,\epsilon)$ is also a k-simplex of $\mathcal{C}_{VR}(X,\epsilon')$, hence $\mathcal{C}_{VR}(X,\epsilon)\subseteq\mathcal{C}_{VR}(X,\epsilon')$. It follows that the following diagram where the vertical maps are inclusions is commutative,
\[
\begin{matrix}
\mathcal{C}_{VR}(X,\epsilon): &\cdots\rightarrow &  C_{n+1} & \xrightarrow{\partial_{n+1}} & C_n & \xrightarrow{\partial_n} & C_{n-1} & \rightarrow\cdots \\
& & \big\downarrow &  & \big\downarrow &  & \big\downarrow & \\
\mathcal{C}_{VR}(X,\epsilon'): &\cdots\rightarrow &  C'_{n+1} & \xrightarrow{\partial'_{n+1}} & C'_n & \xrightarrow{\partial'_n} & C'_{n-1} & \rightarrow\cdots 
\end{matrix}
\]
and hence the inclusion map defines a homomorphism between n-th Homology groups $\mathcal{C}_{VR}(X,\epsilon)$ and n-th homology groups of $\mathcal{C}_{VR}(X,\epsilon')$\cite[Prop 2.9]{Hatcher}. It is then possible for every dimension $n$ to keep track of when (i.e. for what value of $\epsilon$) a new generator is added to the n-th homology group and when it vanishes. This can be depicted in a two dimensional diagram where to each generator corresponds a point $(x,y)$ where $x$ is the value of $\epsilon$ the generator is created and $y$ the value of $\epsilon'$ the generator vanishes. The further away a point is from the main diagonal, the more the corresponding generator survives (see fig.\ref{rips_diagram}). The idea behind Persistent homology is that topological features that survive longer will correspond to actual features of the data and to random noise. For more details one can read articles \cite{ELZ, EH2008, Carlsson2009}. Books \cite{EH2010,Ghrist} cover Persistent Homology as well as a more broad view of applications of Topology. 

Note, the there is another way to depict persistence by using the so-called barcode diagrams\label{barcode_diagrams}. In barcode diagrams, each generator is represented by a horizontal bar with starting point the value of $\epsilon$ the generator is created and end point the value of $\epsilon$ the generator vanishes.The bars are stacked one above the other with the bar that has the smaller starting point at the top (ex. see fig.\ref{two_circles_b})\footnote{Note that the barcode diagrams presented here does not contain bars ending at infinity.}.

\subsection{Mapper}\label{panagopoulos_subsec3_2}
Mapper algorithm was presented by Singh, M\'emoli and Carlsson \cite{SMC}. As the authors state: 
\begin{quote}
The basic idea can be referred to as partial clustering, in
that a key step is to apply standard clustering algorithms to
subsets of the original dataset, and then to understand the
interaction of the partial clusters formed in this way with
each other.
\end{quote}
Mapper algorithm starts with a finite subset $X\subset\mathbb{R}^d$, a map $f:X\rightarrow Z$ to a topological space $Z$ and a covering $\{U_\alpha\}_{\alpha\in A}$ of $Z$\footnote{Usually $Z=\mathbb{R}$ or $Z=\mathbb{R}^2$.}. Then,
\begin{enumerate}\label{mapper_algorithm}
\item every non empty set $f^{-1}(U_\alpha)$ is clustered using a clustering algorithm. Thus, for every $\alpha\in A$ where $f^{-1}(U_\alpha)\neq\emptyset$ we have a finite set $\{C_{\alpha 1},\ldots,C_{\alpha k_\alpha}\}$.
\item a simplicial complex is constructed with vertex set all clusters $\{C_{\alpha j}: 1\leq j\leq k_\alpha, f^{-1}(U_\alpha)\neq\emptyset\}$ and where an n-simplex $\{C_{\alpha_0 i_0},\cdots, C_{\alpha_n i_n}\}$ belongs to the complex if and only if $C_{\alpha_0 i_0}\cap\cdots\cap C_{\alpha_n i_n}\neq\emptyset$.
\end{enumerate}

In the original article of Singh, M\'emoli and Carlsson single-linkage clustering \cite{Wiki_SinglLnk} was used in step (1). More details on Mapper can be found in \cite{SMC,Carlsson2009,Goldfarb}.

\section{Applications of TDA to clustering}\label{panagopoulos_sec4}
\subsection{Using Persistent Homology}\label{panagopoulos_subsec4_1}
It is obvious that Persistent Homology can be used to obtain information about the shape of data. In especial, it can provide information on the number of connected components and help decide whether data have spherical shape or not. It this article we present the following examples:
\begin{itemize}
\item {\bf two squares}: In this example, one hundred points were selected at random using the uniform distribution so that they are in a square with vertices $(0,0),\,(1,0),\,(1,1),\,(0,1)$. And similarly, one hundred points were selected in a square with vertices $(5,5),\,(6,5),\,(6,6),\,(5,6)$. The dataset is depicted in fig.\ref{two_boxes_a}.\\

The corresponding persistent diagram is in fig.\ref{two_boxes_b}. Notice that in dimension 0, there is a group of points near the origin, a point near $(0,6)$\footnote{The coordinates are $(0,t')$, where $t'$ is equal to the closest distance between the two squares which is approximately $4\sqrt{2}$.} and a point with zero x-coordinate on the horizontal dashed line that represents infinity. This is a clear indication that the corresponding Vietoris-Rips complex $\mathcal{C}(X,\epsilon)$ has many connected components for small values of $\epsilon$. Those can be attributed to random noise. The other two points indicate that, at large scale, there exist two connected components that are merged into one for $\epsilon\simeq 6$\footnote{Actually, $\epsilon=t'$.}. For dimension 1, all the points lie close to the diagonal, hence it can be deduced that the dataset does not contain any holes. 

\begin{figure}[ht]
\centerline{
  \subfigure[]
     {\includegraphics[width=1.5in]{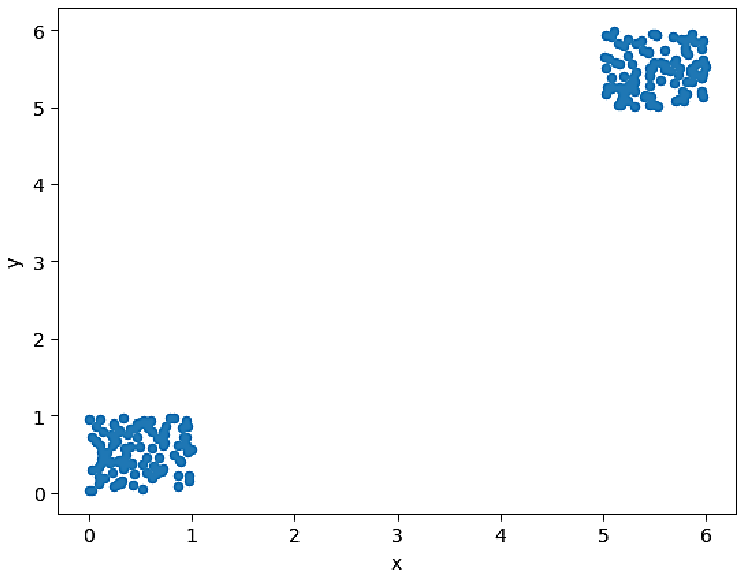}\label{two_boxes_a}}
  \hspace*{4pt}
  \subfigure[]
     {\includegraphics[width=1.5in]{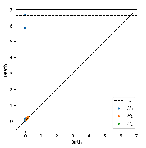}\label{two_boxes_b}}
}
\caption{Two squares made up from random points (a) and the corresponding persistent diagram (b).} \label{two_boxes} 
\end{figure}

\item {\bf two circles}: In this example, five hundred points were selected at random using the uniform distribution so that they are in a ring with inner radius 1 and outer 2. And similarly one thousand points were selected in a ring with inner radius 5 and outer 10. The dataset is depicted in fig.\ref{two_circles}.\\

The corresponding persistent diagram is in fig.\ref{two_circles_a}. In dimension 0, there is a group of points with zero x-coordinate and y-coordinate ranging from zero to approximately one. There is also a point at a height little before two and one lying at the dashed line representing infinity. As in the case of two squares, this is an indication that there are two distinct groups. In dimension 1, we several points above the diagonal that indicate the existence of holes. Most of them are two points one at approximately $(0.4,\,1.8)$ and one at approximately $(0.6,\,3.5)$ (fig.\ref{two_circles_b}) that represent the two circles. This is should act as a warning against using k-means.
\begin{figure}[ht]
\centerline{
  \subfigure[]
     {\includegraphics[width=1.3in]{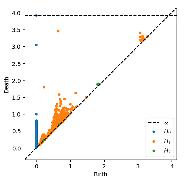}\label{two_circles_a}}
  \hspace*{4pt}
  \subfigure[]
     {\includegraphics[width=1.3in]{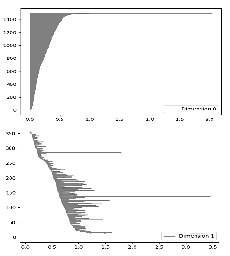}\label{two_circles_b}}
}
\caption{Persistent diagram and barcode diagrams for the two circles dataset.} \label{two_circles_persistent} 
\end{figure}

\item {\bf Iris dataset}: The Iris dataset\cite{Anderson,Fisher} contains measurements of petal and sepal length and width  from one hundred fifty samples from three species (Setosa, Versicolor, Virginica) of iris. There are fifty samples from each of Iris Setosa, Versicolor and Virginica species fig.\ref{iris_pairplot}.
\begin{figure}
\centerline{\includegraphics[width=13cm]{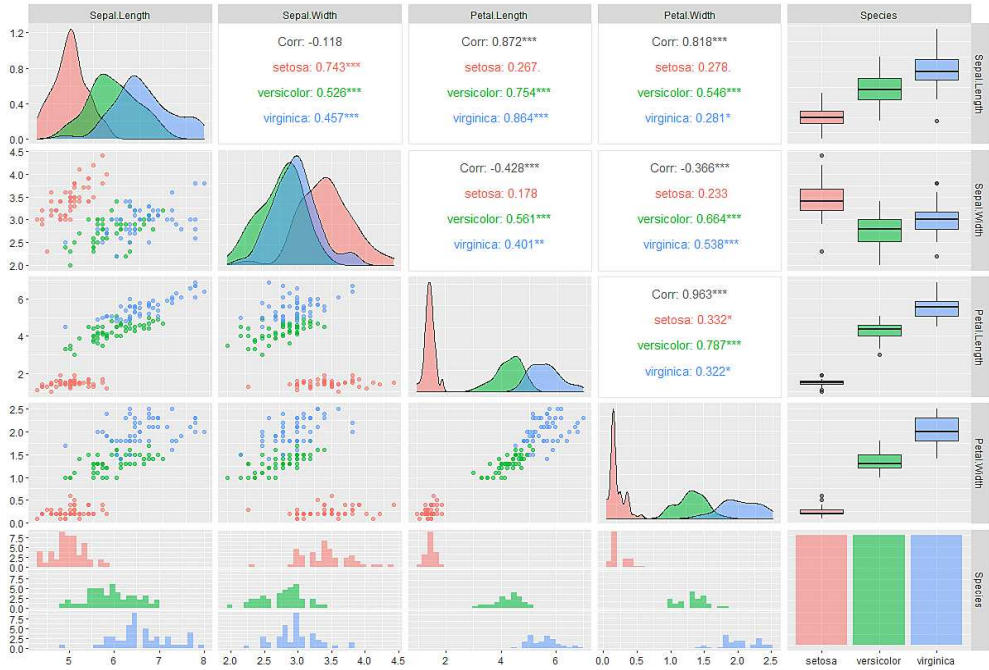}}
\caption{Plot of Iris dataset. In the main diagonal there are the histograms for each variable for each of the three species. Pairwise scatter plots (resp. correlations for the whole dataset and per specie) are below (resp. above) the main diagonal. The last column to the right contains the box-plots for each variable and specie.}
\label{iris_pairplot}
\end{figure}

The corresponding persistent diagram is in fig.\ref{iris_persistent}. In dimension 0, one can clearly see two groups as indicated by the two top points with zero x-coordinate. It is not clear from the diagram if the dataset can be split in one, two or three more groups. The points that correspond to dimensions 1 and 2 lie close to the diagonal, hence there is no indication that k-means will not be an effective clustering method.
\begin{figure}
\centerline{\includegraphics[width=6cm]{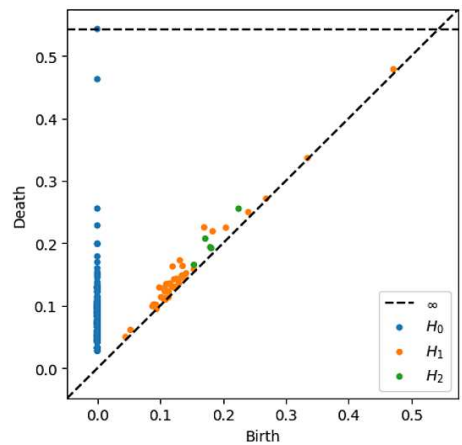}}
\caption{Persistent diagram of the Iris dataset.}
\label{iris_persistent}
\end{figure}

\item {\bf Bank Marketing dataset}: The Bank Marketing dataset\cite{MCR} contains data for approximately forty five thousand customers of a Portuguese bank. Specifically, the data are from a marketing campaign offering term deposits. It contains information on customers age, job, marital status, education, possession of housing or personal loan e.t.c. The standard use of the dataset is for testing classification algorithms. The goal being to classify each customer as either acquiring a term deposit or not. For our purposes we will retain information about   age, job, marital status, education, possession of housing or personal loan and whether there is credit in default. The categorical variables are encoded to ordinal variables (in the case of education) or to groups of binary variables with label or one-hot encoding. After removing cases with null values around thirty thousand records remain. Due to resources restrictions, a random sample of four thousand points is selected for analysis.\\

From the persistence diagram (fig.\ref{bank_marketing_persistent}), it can be seen that there are three or four clusters. Furthermore, it is clear that the dataset contains 1-dimensional holes and thus it is not prone for clustering by k-means.

\begin{figure}
\centerline{\includegraphics[width=6cm]{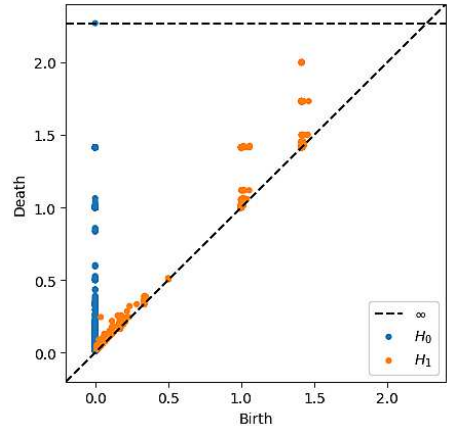}}
\caption{Persistent diagram of the Bank Marketing dataset.}
\label{bank_marketing_persistent}
\end{figure}
\end{itemize}

\subsection{Using Mapper}\label{panagopoulos_subsec4_2}
In this subsection we present the results of application of Mapper algorithm to the above datasets. For the two squares and two circles datasets, the projection in first coordinate was used as map $f:X\rightarrow Z$ (see \ref{mapper_algorithm}). For the Iris dataset and Bank Marketing Datasets projection on principal components derived by PCA was used. It should be noted that since this is an introduction, we will not present an exhaustive analysis of the results.
\begin{itemize}
\item {\bf two squares}: In fig.\ref{mapper_two_squares} we can see the result of Mapper on the dataset where the points are selected at random within two squares. Each square is represented by two connected nodes.
\begin{figure}
\centerline{\includegraphics[width=8cm]{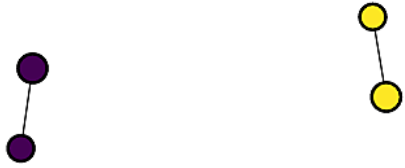}}
\caption{Application of Mapper algorithm on two squares dataset.}
\label{mapper_two_squares}
\end{figure}

\item {\bf two circles}: While in the case of two squares Mapper represents the two separate clusters of points, this is not the case for the two circles. In fig.\ref{mapper_two_circles} we can see the results when using DBSCAN or k-means for clustering. The expected result would be a graph consisting of nodes forming two cycles. This was not possible to achieve even thought a variety of parameters and clustering methods were tested. Figure \ref{mapper_two_circles} is an indication of the high sensitivity of the Mapper algorithm to the selection of clustering algorithms.  
\begin{figure}
\centerline{
  \subfigure[]
     {\includegraphics[width=2.5in]{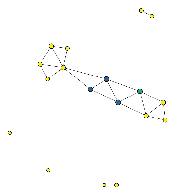}\label{mapper_two_circles_a}}
  \hspace*{4pt}
  \subfigure[]
     {\includegraphics[width=2.5in]{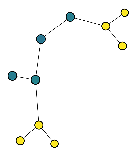}\label{mapper_two_circles_b}}
}
\caption{Application of Mapper to the two circles. Projection to the first coordinate was used. For the image to the left DBSCAN was used as clustering algorithm, while k-means was used for the image to the right. Colour indicates ratio between points from the outer and points from the inner circles in each node. Ranging from yellow for nodes where all points are from the outer circle to blue for nodes with all points from the inner circle.} \label{mapper_two_circles} 
\end{figure}

\item {\bf Iris dataset}: For Iris dataset, projection on the first two principal components was used as map $f:X\rightarrow Z$ (\ref{mapper_algorithm}).It is clear from the the pair plots that, Setosa, one of the three species of Iris is clearly separated from the other two (see fig.\ref{iris_pairplot}). Mapper manages to capture this fact. As we can see Mapper constructs a graph with two connected components. One of them corresponds to Setosa. The other corresponds to the other two species. Furthermore, in the second connected component nodes in one side correspond to Versicolor and nodes on the other side to Virginica. 
\begin{figure}
\centerline{\includegraphics[width=8cm]{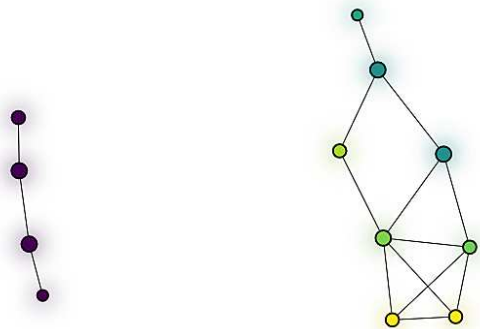}}
\caption{Application of Mapper to Iris dataset. Colour indicates species. Purple for Setose, green for Versicolor and yellow for Virginica.}
\label{iris_mapper}
\end{figure}

\item {\bf Bank Marketing dataset}: 
For Bank Marketing dataset, projection on the first five principal components was used as map $f:X\rightarrow Z$ (see \ref{mapper_algorithm}). Mapper reveals, several well defined groups of customers (see fig.\ref{bank_marketing_mapper}). By examining the values of the variables we are able to determine which of them are used to define those groups. This process reveals that the existence of a housing loan or the marital status play a major role. In contrast age does not. 
\begin{figure}
\centerline{\includegraphics[width=8cm]{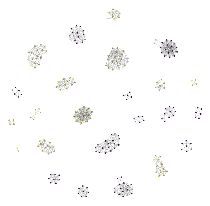}}
\caption{Application of Mapper to Bank Marketing dataset.}
\label{bank_marketing_mapper}
\end{figure}
\end{itemize}
\section{Concluding remarks}\label{panagopoulos_sec5}
In the previous sections we applied Persistent Homology and Mapper to some simple datasets. In particular, he experiments we conducted indicate that Persistent Homology can be used successfully to study the number of clusters and the shape of a data set. In contrast, Mapper algorithm produced mixed results. In some cases managing to capture key information about the data while in some other cases filing to do so.   

More applications, can be found in \cite{Carlsson2009,EH2010,ELZ}. In \cite{GKS} the authors use Topological Data Analysis to detect relationships between products sold on a local level and products sold on a national level. It is a very interesting idea that it also highlights on of the obstacles in using Topological Data Analysis, the fact that the current algorithms do not scale to accommodate for Big Data. In \cite{Lumetal} the authors apply Mapper to analyse gene expression in breast tumours, performance data from the NBA and voting data from US House of Representatives. In \cite{KMJ2017} the authors compare Mapper against k-means and hierarchical clustering on the the task of image popularity in social media.

It is clear from the applications above, that Persistence Homology and the Mapper algorithm can be used to provide insights on datasets. Unfortunately, there are a couple of reasons that make wider adoption of those tools difficult. Some of them are the rather technical background required to understand those tools, the existence of some well established techniques for addressing similar problems and the lack of high efficient algorithms that scale well to Big Data. That said, the field is active with research going on both on applications and on the theoretical background (ex. see \cite{Casas,PP,BP}).

%\bibliographystyle{ws-rv-van}

%\printindex[aindx]                 % to print author index
%\printindex                         % to print subject index

\begin{thebibliography}{9}
\bibitem{Akhanli}
   S. E. Akhani, \emph{Distance construction and clustering of
football player performance data},
   PhD Thesis, (2019) \url{https://discovery.ucl.ac.uk/id/eprint/10065964/1/thesis.pdf}. 
  
\bibitem{ADHP}
	D. Aloise, A. Deshpande, P. Hansen, P. Popat, \emph{NP-hardness of Euclidean sum-of-squares clustering}, \emph{Machine Learning}, Springer, {\bf 75} (2), (2009), pp. 245--248.
	
\bibitem{Anderson}
	E. Anderson, \emph{The irises of the Gaspe Peninsula}, \emph{Bulletin of the American Iris Society}, {\bf 59}, pp. 2--5 (1935).

\bibitem{AV}
   D. Arthur, S. Vassilvitskii, \emph{k-means++: The advantages of careful seeding}, \emph{Proceedings of the eighteenth annual ACM-SIAM symposium on Discrete algorithms}, Society for Industrial and Applied Mathematics, (2007), \url{http://ilpubs.stanford.edu:8090/778/1/2006-13.pdf}. 
   
\bibitem{BP}
   S. Basu, L. Parida \emph{Spectral sequences, exact couples and persistent homology of filtrations}, \emph{Expositiones Mathimaticae}, {\bf 35} 1, pp. 119--132, (2017), \url{https://www.sciencedirect.com/science/article/pii/S0723086916300378}.

\bibitem{Carlsson2009}
   G. Carlsson, \emph{Topology and Data}, \emph{Bull. Amer. Math. Soc.}, {\bf 46}, pp. 255--308, (2009).
   
\bibitem{Casas}
   A.T. Casas, \emph{Distributing Persistent Homology via Spectral Sequences},
   \emph{arXiv}, (20202), \url{https://arxiv.org/abs/1907.05228}.

\bibitem{EH2008}
	H. Edelsbrunner, J. Harer, \emph{Persistent Homology - a survey}, \emph{Contemporary Mathematics}, AMS, {\bf 453}, pp. 257--282, (2008), \url{https://www.maths.ed.ac.uk/~v1ranick/papers/edelhare.pdf}.
	
\bibitem{EH2010}
	H. Edelsbrunner, J. Harer, \emph{Computational Topology An introduction}, \emph{AMS}, (2010), \url{https://www.maths.ed.ac.uk/~v1ranick/papers/edelcomp.pdf}.
	
\bibitem{ELZ}
	H. Edelsbrunner, D. Letscher, A. Zomorodian, \emph{Topological Persistence and Simplification}, \emph{Discrete Comput. Geom.}, {\bf 28}, pp. 511--533, (2002), \url{https://link.springer.com/content/pdf/10.1007/s00454-002-2885-2.pdf}.

\bibitem{Fisher}
	R. A. Fisher, \emph{The use of multiple measurements in taxonomic problems}, \emph{Annals of Eugenics}, {\bf 7}, Part II, pp. 179--188 (1936), \url{https://archive.ics.uci.edu/ml/datasets/iris}.

\bibitem{Ghrist}
	R. Ghirst, \emph{Elementary Applied Topology}, \emph{Self publication}, (2014), \url{https://www2.math.upenn.edu/~ghrist/notes.html}.

\bibitem{GKS}
   A. Goldfarb, J. B. Kwon, T. Snider, \emph{Detecting potential product segments using topological data analysis}, \emph{Working paper}, (2016) \url{https://www-2.rotman.utoronto.ca/~agoldfarb/TDA.pdf}.
   
\bibitem{Goldfarb}
   B. Goldfarb, \emph{The Mapper algorithm and its applications},
   \emph{15th Annual Workshop on Topology and Dynamical Systems}, (2018), \url{http://topology.nipissingu.ca/workshop2018/slides/Goldfarb-Data-Beamer.pdf}.
   
\bibitem{Hatcher}
	A. Hatcher, \emph{Algebraic Topology}, \emph{Cambridge University Press}, (2002), \url{https://pi.math.cornell.edu/~hatcher/AT/ATpage.html}.

\bibitem{KR}
	L. Kaufman,  P. J. Rousseeuw, \emph{Partitioning Around Medoids (Program PAM)}, \emph{Wiley Series in Probability and Statistics}, Hoboken, NJ, USA: John Wiley \& Sons, Inc., (1990), pp. 68--125.
	
\bibitem{KMJ2017}
   A. Khaled, K. Minkyu, L. Jeongkyu, \emph{Extracting Knowledge from the Geometric Shape of Social Network Data Using Topological Data Analysis},
   \emph{Entropy}, {\bf 19} (7), (2017), \url{https://www.mdpi.com/1099-4300/19/7/360/htm}. 
   
\bibitem{Kleinberg}
   J. M. Kleinberg, \emph{An impossibility theorem for clustering},
   \emph{NIPS}, pp. 446--453, (2002), \url{https://www.cs.cornell.edu/home/kleinber/nips15.pdf}. 

\bibitem{Lumetal}
	P. Y. Lum et al., \emph{Extracting insights from the shape of complex data using topology},  \emph{Scientific Reports}, Nature, {\bf 3} 1236, (2013), \url{https://www.nature.com/articles/srep01236}.

\bibitem{MNV}
	M. Mahajan, P. Nimbhorkar, K. Varadarajan, \emph{The Planar k-Means Problem is NP-Hard},  \emph{Lecture Notes in Computer Science}, Springer, {\bf 5431}, (2009), pp. 274--285.

\bibitem{MCR}
	S. Moro, P. Cortez, P. Rita, \emph{A Data-Driven Approach to Predict the Success of Bank Telemarketing}, \emph{Decision Support Systems}, Elsevier, {\bf 62}, pp. 22--31, (2014) \url{https://archive.ics.uci.edu/ml/datasets/bank+marketing}.

\bibitem{NLC}
	M. Nicolau, A. J. Levine, G. Carlsson, \emph{Topology based data analysis identifies a subgroup of breast cancers with a unique mutational profile and excellent survival}, \emph{Proceedings of the National Academy of Sciences of the United States of America}, {\bf 108} 17, pp. 7265--7270, (2011), \url{https://www.ncbi.nlm.nih.gov/pmc/articles/PMC3084136/}.
	
\bibitem{ORSS}
   R. Ostrovsky, Y. Rabani, L. Schulman, C. Swamy, \emph{The Effectiveness of Lloyd-Type Methods for the k-Means Problem}, \emph{Proceedings of the 47th Annual IEEE Symposium on Foundations of Computer Science (FOCS'06)}, IEEE, (2006), pp. 165–174, \url{https://web.cs.ucla.edu/~rafail/PUBLIC/76.pdf}.

\bibitem{Pedregosa}
   F. Pedregosa et al., \emph{Scikit-learn: Machine Learning in Python},
   \emph{JMLR}, {\bf 12} (85), (2011), \url{https://jmlr.csail.mit.edu/papers/v12/pedregosa11a.html}.

\bibitem{PP}
	M. Piekenbrock, J. A. Perea, \emph{Move Schedules:Fast persistence computations in sparse dynamic settings}, \emph{arXiv}, (2021), \url{https://arxiv.org/abs/2104.12285}

\bibitem{ST}
	N. Saul, C. Tralie, \emph{Scikit-TDA: Topological Data Analysis for Python}, (2019), \url{https://github.com/scikit-tda/scikit-tda}

\bibitem{Sculley}
	D. Sculley, \emph{Web Scale K-Means clustering}, \emph{Proceedings of the 19th international conference on World wide web}, (2010), \url{https://www.eecs.tufts.edu/~dsculley/papers/fastkmeans.pdf}

\bibitem{SMC}
	G. Singh, F. Mémoli, G. Carlsson, \emph{topological Methods for the Analysis of High Dimensional Data Sets and 3D Object Recognition}, \emph{Eurographics Symposium on Point-Based Graphics}, (2007), \url{https://research.math.osu.edu/tgda/mapperPBG.pdf}.

\bibitem{Sot2003}
   Ch. Sotiriou et al., \emph{Breast cancer classification and prognosis
based on gene expression profiles from a population-based study},
   \emph{PNAS}, {\bf 100} (18), (2003), \url{https://www.pnas.org/content/pnas/100/18/10393.full.pdf}.

\bibitem{YMJ}
   L. Yuanhong, D. Ming, H. Jing, \emph{A Gausssian Mixture Model to Detect Clusters Embedded in Feature Subspace},
   \emph{Commun. Inf. Syst.}, {\bf 7} (4), pp. 337--352, (2007), \url{https://projecteuclid.org/journals/communications-in-information-and-systems/volume-7/issue-4/A-Gaussian-Mixture-Model-to-Detect-Clusters-Embedded-in-Feature/cis/1211574970.full}.

\bibitem{WDM}
   H. Wagner, P. Dlotko, M. Mrozek \emph{Computational topology in text mining},
   \emph{Computational Topology in Text Mining. In: Ferri M., Frosini P., Landi C., Cerri A., Di Fabio B. (eds) Computational Topology in Image Context}, Lecture Notes in Computer Science, {\bf 7309}, Springer, Berlin, Heidelberg, (2012),  \url{https://www2.math.upenn.edu/~dlotko/textMining.pdf}.

\bibitem{Wiki}
   Wikipedia, \emph{Cluster Analysis},
   \emph{Wikipedia}, \url{https://en.wikipedia.org/wiki/Cluster_analysis}.
   
\bibitem{Wiki_SinglLnk}
   Wikipedia, \emph{Single-linkage clustering},
   \emph{Wikipedia}, \url{https://en.wikipedia.org/wiki/Single-linkage_clustering}.
\end{thebibliography}
\end{document}